# Following Knots Down Their Energy Gradients

by Louis H. Kauffman
Math, UIC
851 South Morgan Street
Chicago, Illinois 60607-7045
<kauffman@uic.edu>

**Abstract:** This paper details a series of experiments in searching for minimal energy configurations for knots and links using the computer program KnotPlot.

## I. Introduction

This paper details remarkable experiments in simulated physical movement of knots under self-repelling forces. We use the program KnotPlot of Robert Scharein [4,5], making it possible for readers to try these experiments themselves. The experiments are explained in the sections below and they show computer-empirical instances of local energy minima that are not global energy minima.
These examples are obtained by starting with knots or links in particular geometric configurations and letting them evolve and self-repel using the computer program. The paper is organized as follows. Section II gives general background information for the experiments. Section III describes an experiment with the trefoil knot where the starting configurations for the (2,3) and (3,2) versions of the knot lead to different local minima in the experiment. In Section IV we do a similar experiment with a torus link. In Section V we discuss ropelength and examine the behaviour of the ropelength of the minimal energy configurations. Finally in Section VI we show are very remarkable example with the (3,4) torus knot where the undamped and unperturbed self-repulsion leads to a local minimum that is highly symmetrical and not the global minimum. A perturbation of this configuration leads to what is apparently the global minimum. The section ends with discussion and questions. It gives the author pleasure to thank Slavik Jablan for the invitation to write this paper.

## II. Introduction to the Experiment.

Recall that a knot is an embedding of a circle in three dimensional space. The topological type of the knot is its equivalence class under ambient isotopy where two knots are said to be ambient isotopic if one can place them in a continuously varying family of knots starting with one, and ending with the other. We say that a knot is unknotted if its ambient isotopy class contains a simple planar circle, the unknot.

In general it is a difficult problem to tell whether a knot is knotted or unknotted. There are combinatorial algorithms that can determine knottedness from a graphical representation of a knot. On the other hand, there are physically motivated algorithms that apparently can unknot a knot. One of these is motivated by the following idea: Coat the knot with electrical charge and then let it self-repel without changing its length. If it is unknotted, one hopes that the self-repulsion will push it apart from itself into an obviously unknotted form.

This paper recounts experiments that the author has performed on self-repelling knots using the computer program KnotPlot [4, 5]. KnotPlot is a program written by Rob Scharein and it models knots that self-repel with a $(1/r)^d$ force where the exponent d is adjustable between 2 and 6. Here r stands for the radial distance between two points on the embedded curve that represents the knot.

In the computer model, the knot is represented by a string of vectors $(v_1, v_2, ...., v_n)$ so that the embedding for the knot is given by connecting the tips these vectors cyclically so that $v_1$ is connected to $v_2$, $v_2$ to $v_3$ and so on with $v_n$ connected to $v_1$. Thus the knot itself contains the points corresponding to the vectors $v_k$. To implement a self-repelling force on the knot, one computes the forces between $v_k$ and $v_l$ for all distinct pairs k and l in the set {1,2,..., n} and uses these forces to move the points by a small amount that does not cause any of the connecting segments to cross through one another. This process is iterated, and the knot moves in a simulation of self-repulsion.

On can also calculate the Simon energy E(K) [1]for a given knot embedding. This energy corresponds (by definition) to the energy potential for a $(1/r)^2$ force (d=2). E(K) is obtained by summing $1/|v_i - v_j|$ for all distinct pairs i and j in {1,2,...,n}. The program KnotPlot can be asked to display this energy. As the knot self-repels, this energy function tends to a minimum but under certain circumstances will oscillate or change quite slowly.

Part of the circumstances related to the behaviour of the energy functional have to do with an additional feature of the program KnotPlot. Along with implementing a repelling force, the program also models the connections between the points $v_k$ and $v_{k+1}$ as springs with restoring force behaving according to Hooke's law (force proportional to the extension of the length of the distance beyond a base distance corresponding to the length of the contracted spring). One can operate the program so the the springs are either damped or undamped. With damped springs the program operates as though the springs are fully contracted, and no energy is exchanged with the springs. In undamped mode, the knot exchanges energy with its springs, moves around, oscillates and the energy functional will be seen to go up and down. Putting a knot in this process in the undamped mode is a way to subject it to perturbation that may allow it more freedom of movement than it can have in the damped mode. We will see how this plays out in the experiments.

### III. The Trefoil Experiment
There are two versions of the trefoil knot when you regard it as a torus knot. You view it as a (2,3) torus knot, winding twice around the torus in the longitudinal direction and three times around in the meridianal direction, as depicted in Figure 1.

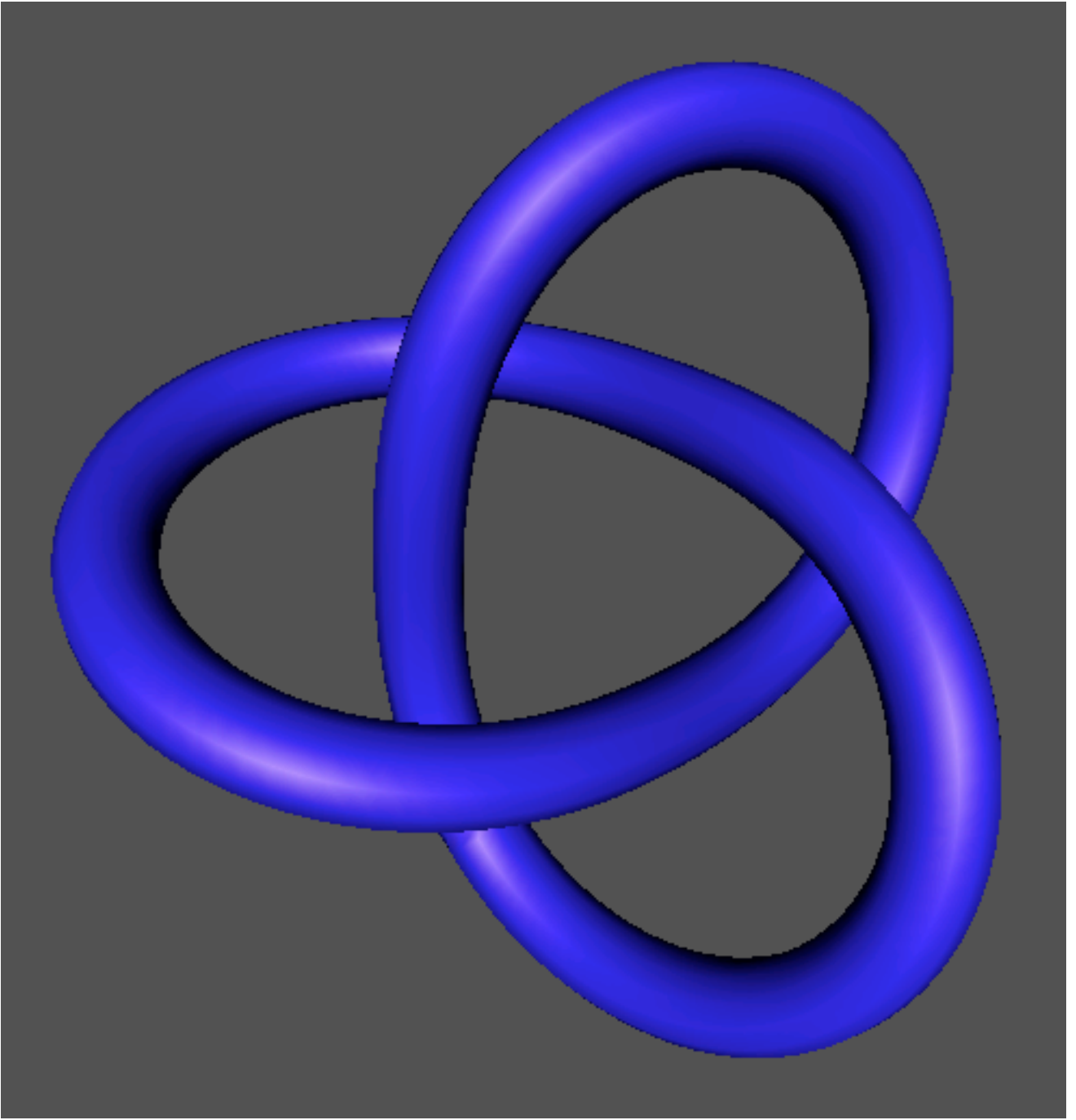

**Figure 1 - The (2, 3) torus knot with E = 169.04011473**

In Figure 1, we have illustrated the image of the torus knot (2,3) that results from running damped self-repulsion until the energy stabilizes at a minimal value. There are 80 vector points in this version of the knot.

The other version of the trefoil knot is the (3,2) torus knot as shown in Figure 2. The (3,2) torus knot winds around three times in the longitude direction for two times in the meridian direction. The (2,3) and (3,2) torus knots are ambient isotopic as knots in three dimensional space. They represent different starting points for the self-repulsion program. Our experiment will proceed from the versions of the (3,2) knot shown in Figure 2. It is at a much higher energy than its eventual minimal configuration and we shall follow it toward that minimal configuration. We will find that this minimal configuration is still at a higher energy than the form shown in Figure 1.

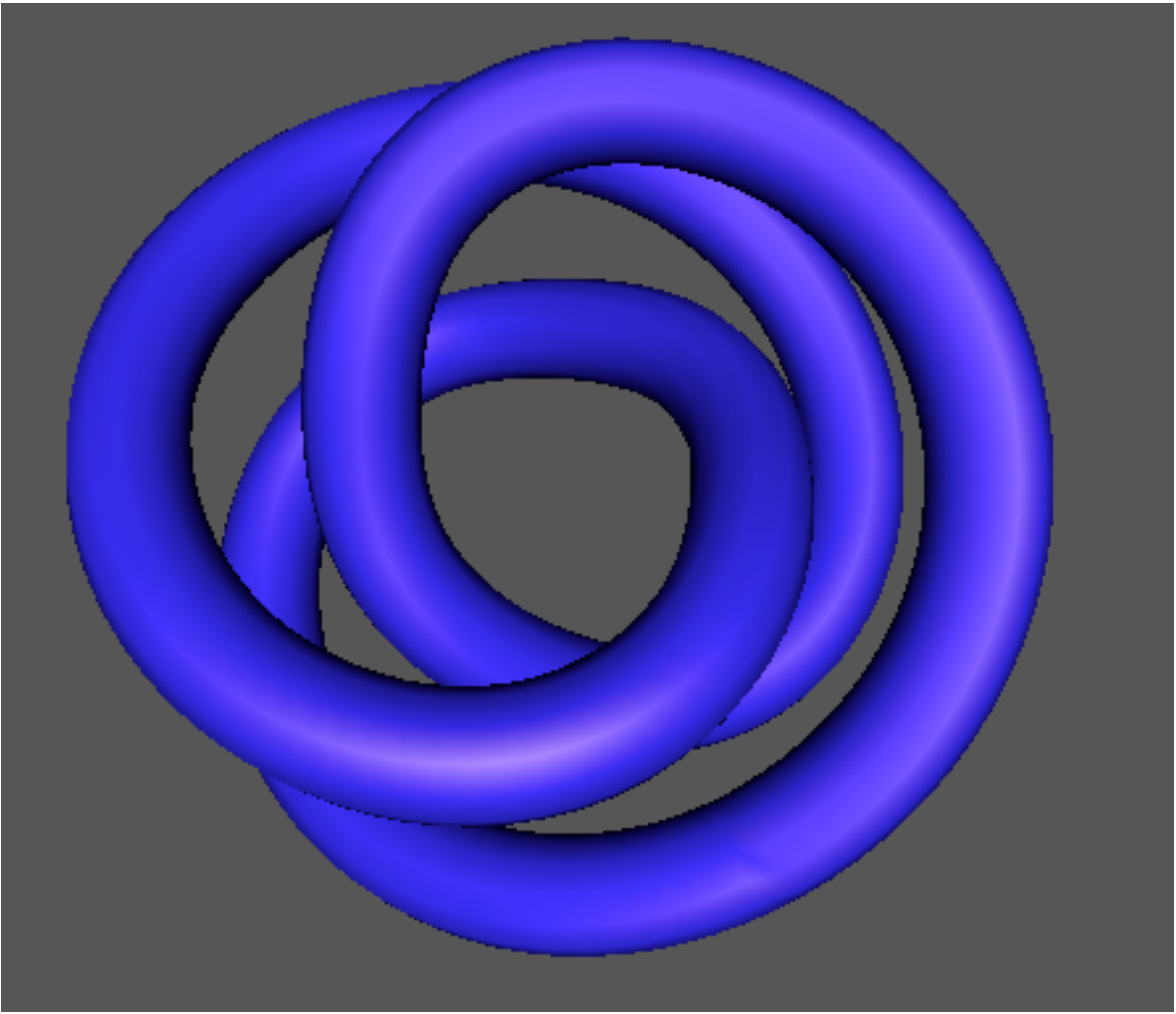

**Figure 2 - The (3,2) torus knot at E= 212.49077.**

Figure 3 shows KnotPlot's minimal energy form using damped self-repulsion for the (3,2) torus knot.

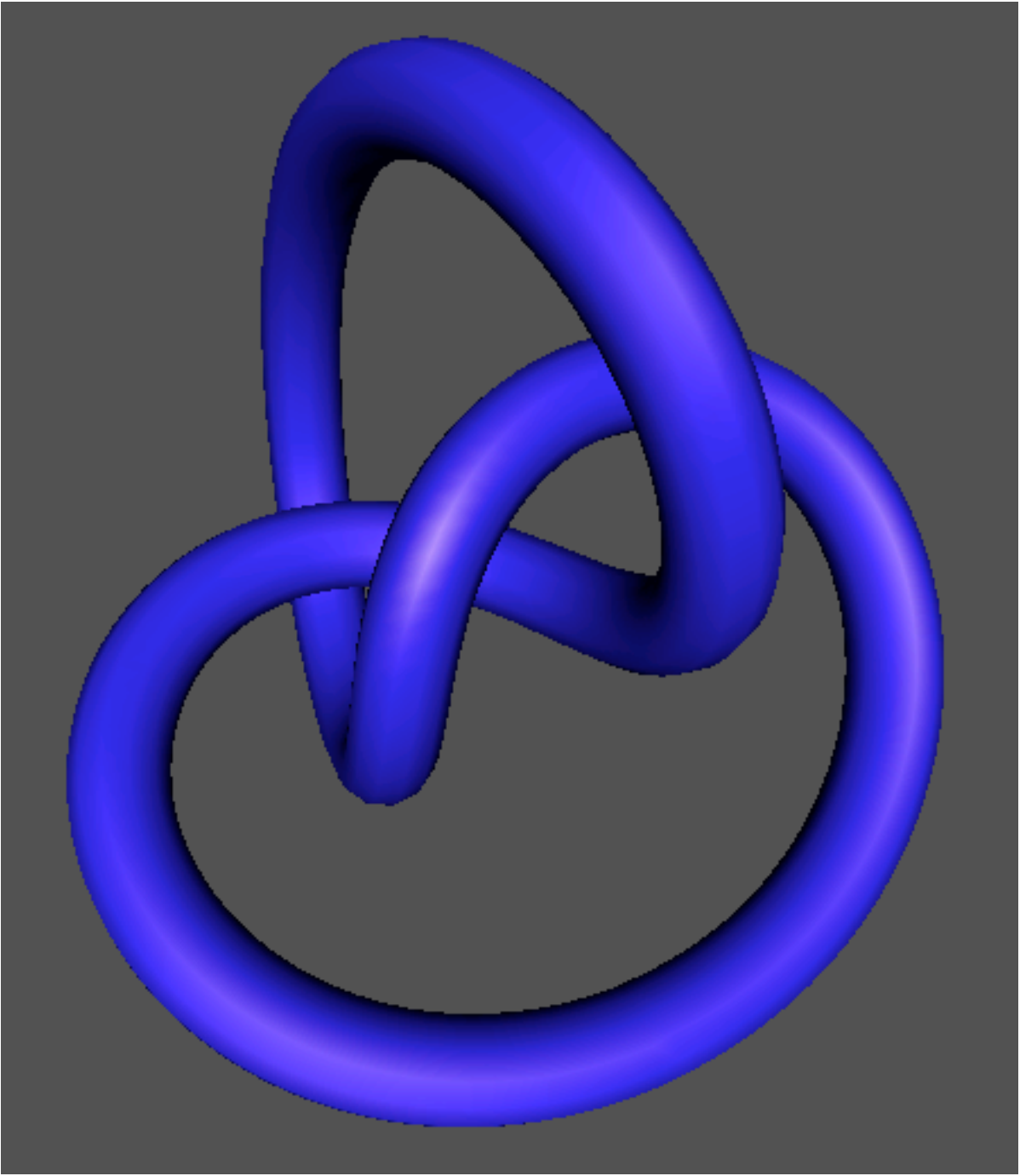

**Figure 3 - The end result for the (3,2) torus knot.
E = 175.074081.**

In Figure 3 we see the result of running the self-repulsion program for until the energy E has stabilized at a specific value and the geometric evolution of the knot appears to be stopped. Figures 4, 5,

6 and 7 show the further evolution of this configuration when we use undamped self-repulsion. The knot then exchanges energy with the springs and goes into an unstable oscillation that results in a descent to a lower energy level. The final result in Figure 7 is the same as the level in Figure 1.

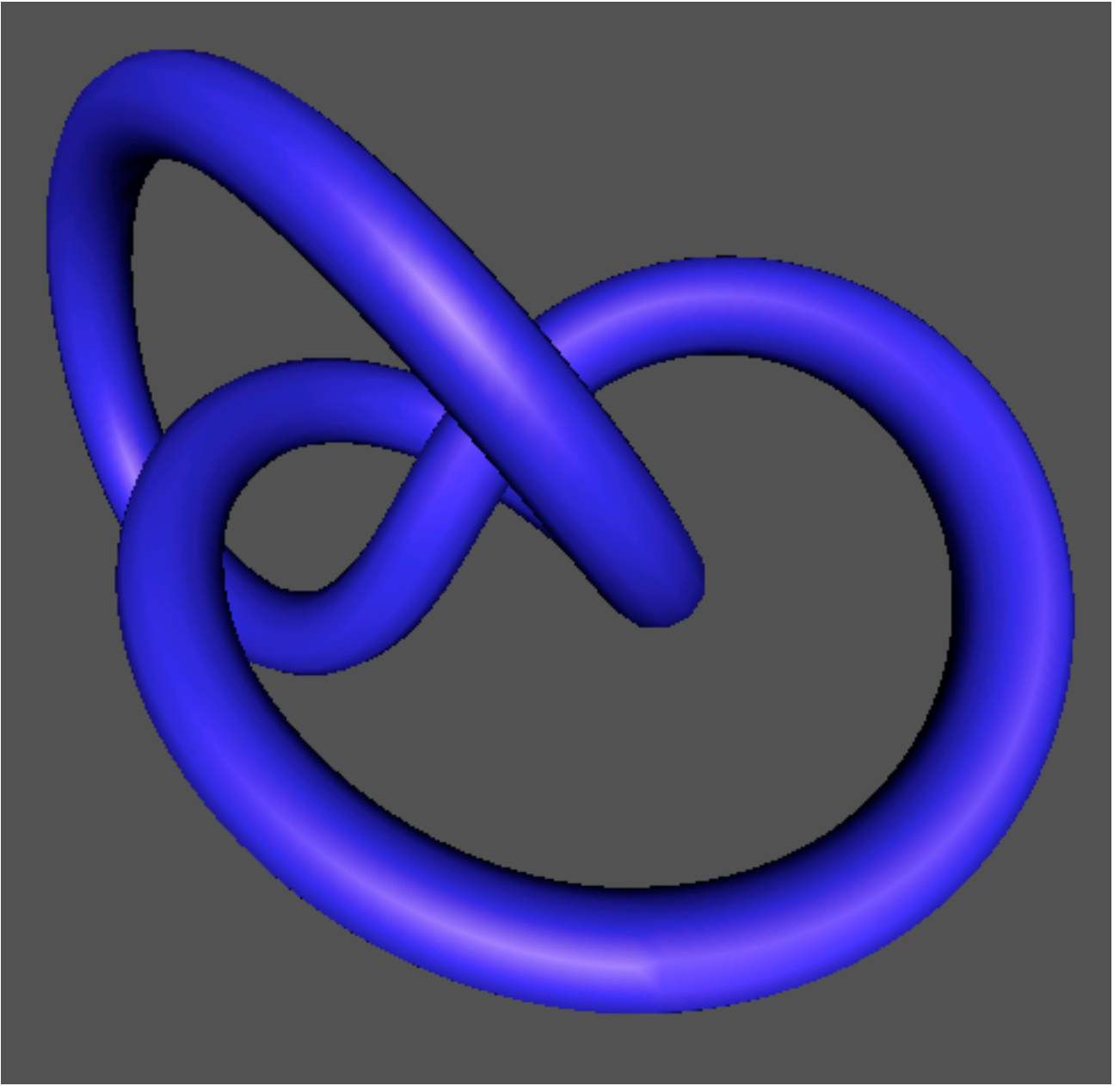

**Figure 4 - E = 173. 034119.**

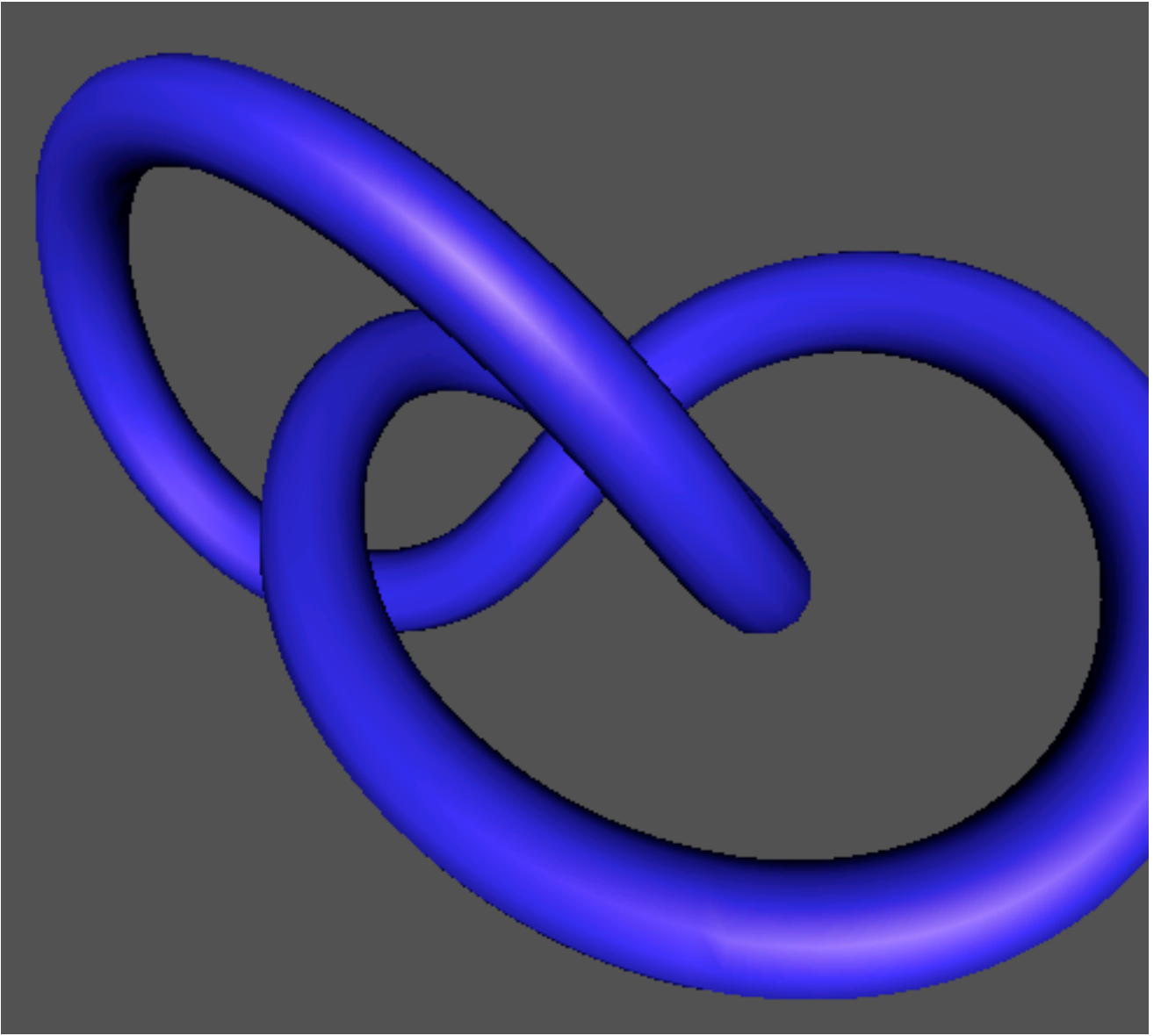

**Figure 5 - E = 171.196564**

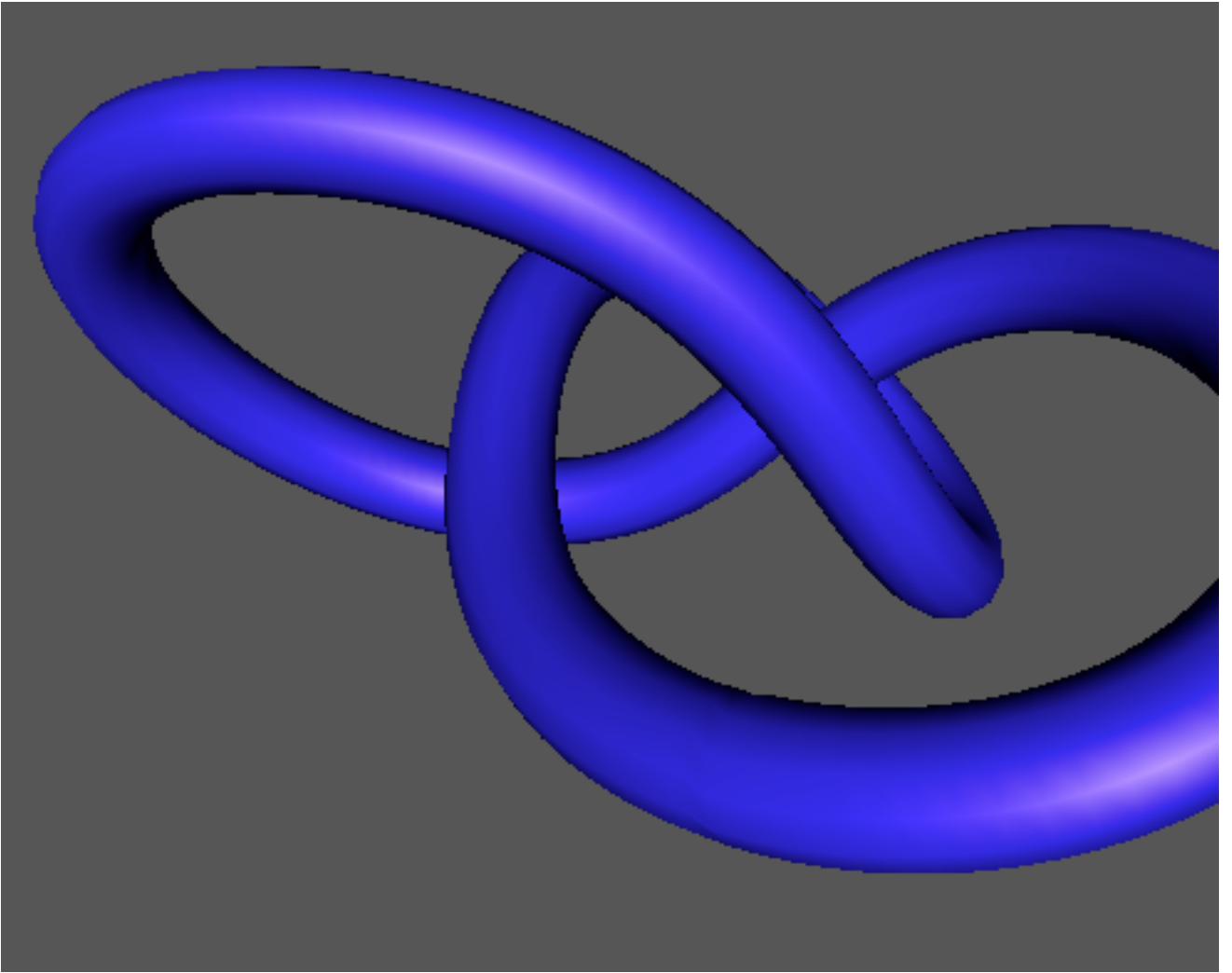

**Figure 6 - E = 169.545120**

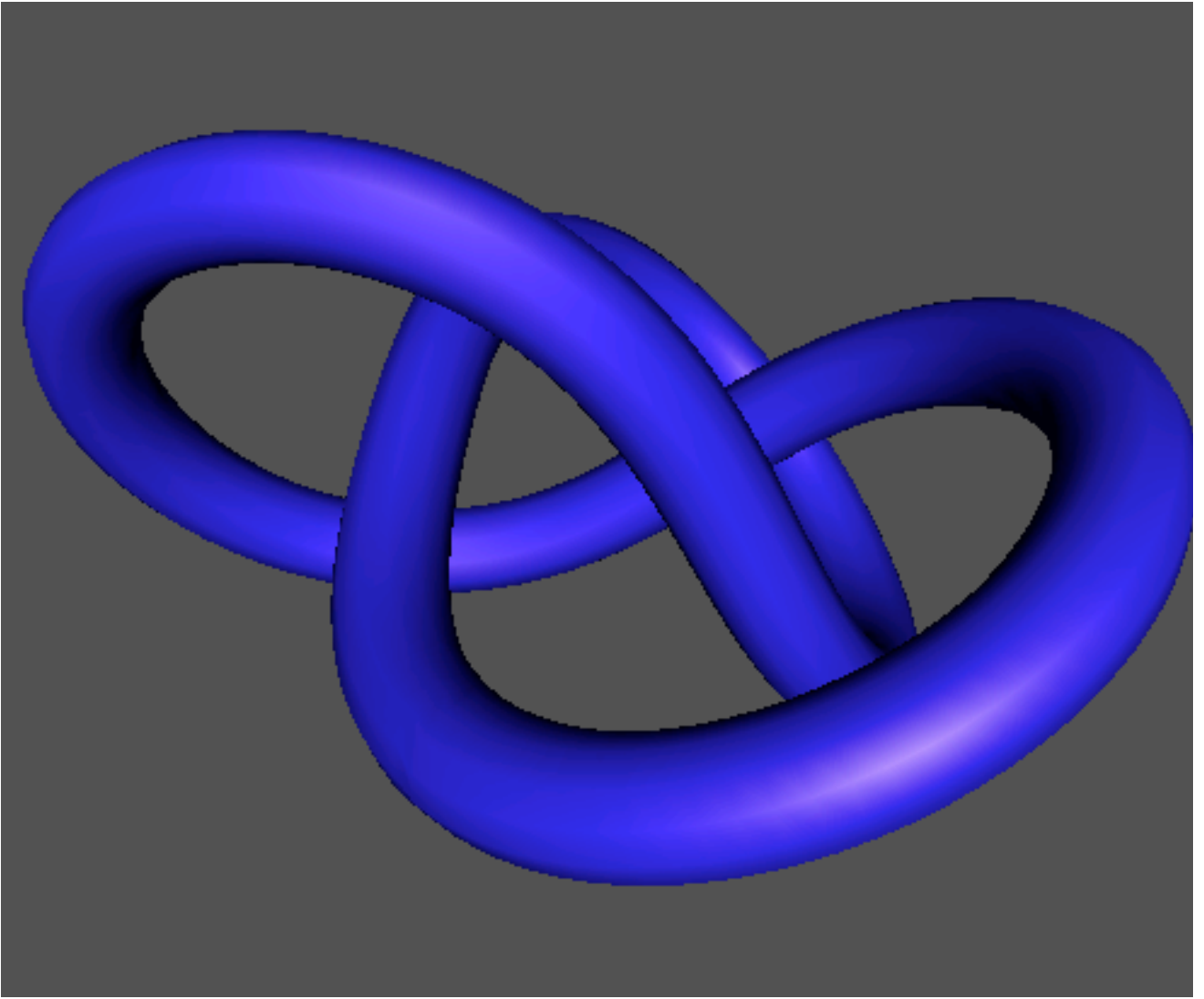

**Figure 7 - E = 169.041534**

Note that this final result in Figure 7 has the same energy as the configuration in Figure 1. In fact if rotated it would be identical with Figure 1.

The apparently stable intermediate stage obtained from the (3,2) torus link as shown in Figure 3 may be indicate a local minimum for the energy for the trefoil knot or this configuration could be the analog of an inflection point, where the rate of evolution becomes very slow. What we can say is that it is for this particular computer program an apparent local minimum, and it gives us a good graphic demonstration of how a stable energy level may, under perturbation, fall into an even lower energy level through further evolution of the form.

## IV. The Torus Link Experiment

For the (4,2) torus link, we find a very similar phenomenon. Starting in the (4,2) configuration the link descends to the form shown in Figure 8.

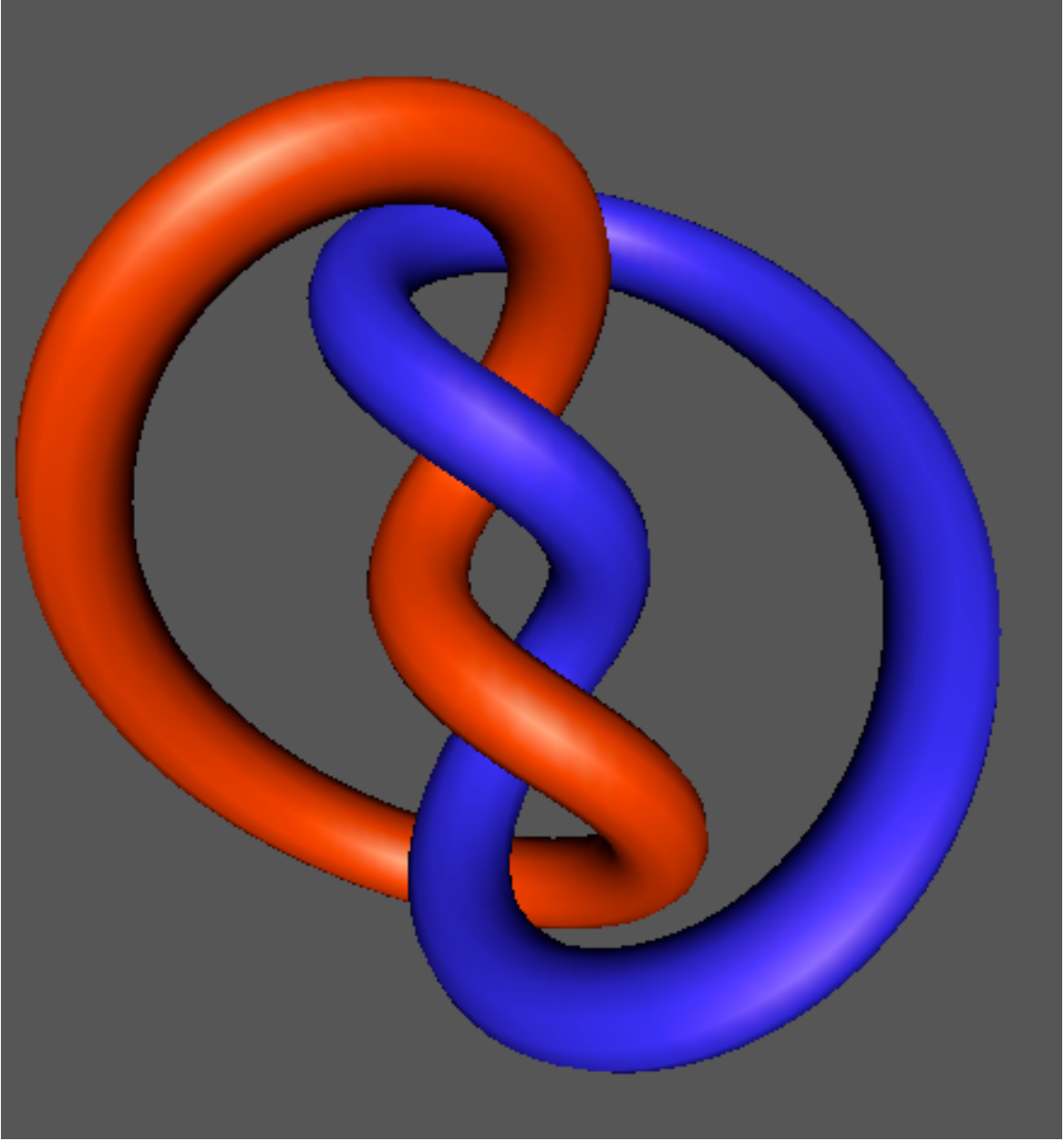

**Figure 8 - (4,2) torus link with E = 199.353058.**

And then on perturbation through the undamped repelling process, this relatively stable form goes into oscillation and descends to the form shown in Figure 10. The form in Figure 9 is the same as the form that results from self-repelling the (2,4) torus link.

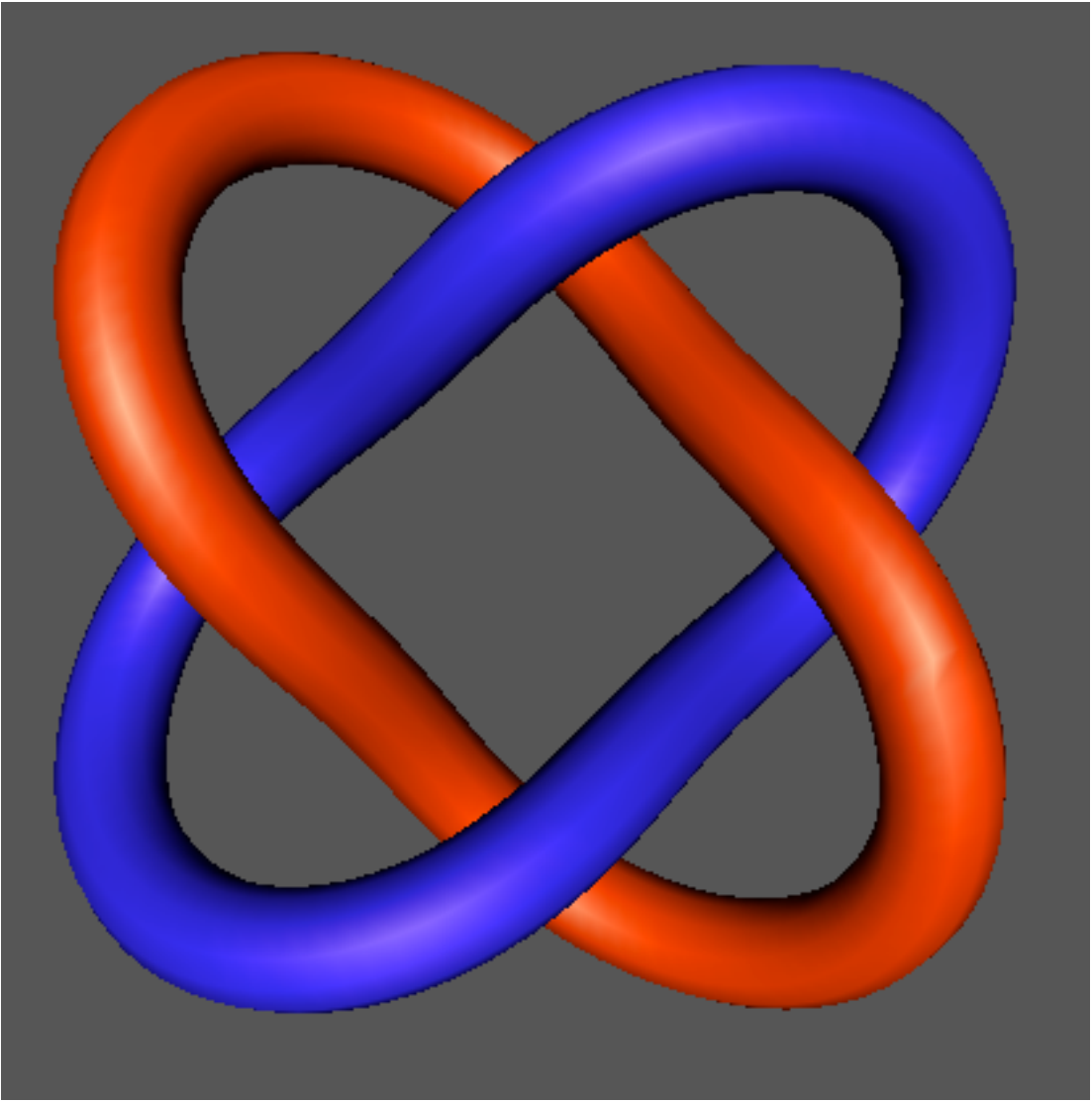

**Figure 9 - Full Descent of the (4,2) Torus Link, E = 185.343277.**

## V. Electrical RopeLength

We now turn to a way to geometrically measure configurations of knots in three-dimensional space that is related to the concept of *minimum rope-length*. The minimum rope length MinRL(K) of a knot K is the minimum over all tubular embeddings of K in three-space of the ratio of the length of the tube to its cylindrical radius. For a specific embedding of a knot, we can increase the cylinder radius until the tube just touches itself. The ratio of tubular length to cylinder radius for such a cylinder is defined to be the *ropelength for this particular embedding*. See [2,3] for more information about these definitions. As in [3] programs have been written to approximate the minimal ropelength for a knot by simulating the process that one can perform on a rope of pulling the rope tight in order to concentrate the knotting on the least amount of rope. The natural question in terms of rope is: What is the least length of rope, for a given diameter, needed to make a given knot-type?

Here we point out that since there is the ropelength RL(K) (defined as above by taking the maximum cylinder radius for the embedding) for any specific embedding K in three-space, we can consider RL(K) as a measure of the complexity of an embedding and ask about the value of RL(K) for an embedding of K that minimizes the knot energy or that is an apparent minimum in a self-repelling process such as considered in this paper. Thus we shall define ERL(K, d) to be the ropelength of a configuration for K that is extremal for the self-repelling process with a force field proportional to $r^{-d}$. We call this the *electrical ropelength* of the knot K. The electrical ropelength is an interesting parameter to study and it may, only for large d, be approximately equal to the standard ropelength. The program we are using calculates the tubular length and the cylinder radius and allows variation of the cylinder radius. Thus we can experiment with the values of the electrical ropelength and the geometry of the resulting tubular knots. In Figure 11 below we illustrate an approximation for the electrical ropelength ERL(T,6) where T is a trefoil knot and the force exponent is equal to 6.

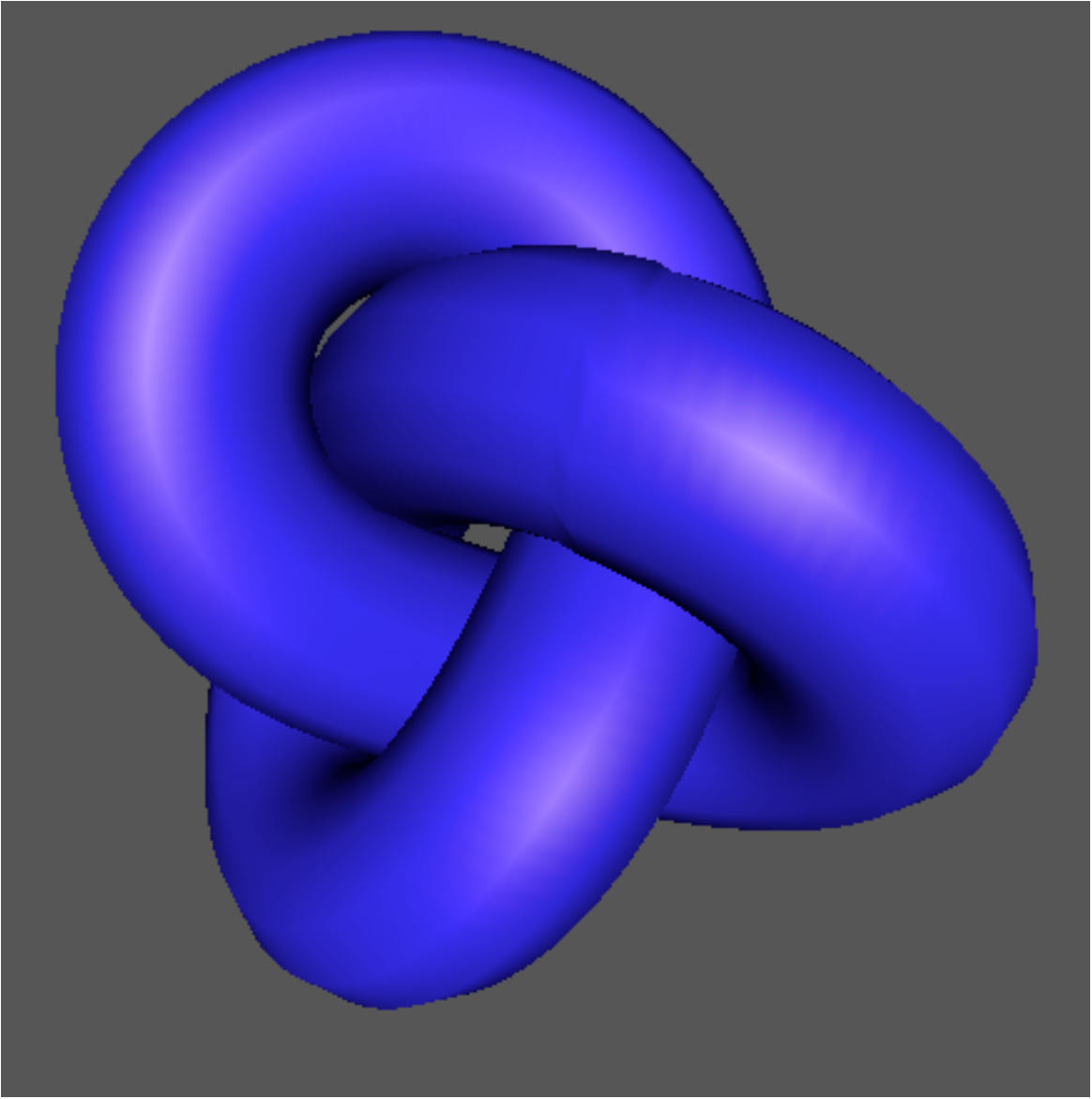

Figure 10- ELR(T,6) ~ 32.68. Electrical Ropelength for the Trefoil Knot.

In practice d=6 seems to be a good exponent both for force repelling experiments and for the electrical ropelength. For smaller d once finds that the ropelength increases. For example, see Figure 11 where we illustrate the trefoil with d = 3.5. There we have ELR(T,3.5) ~ 40.45.

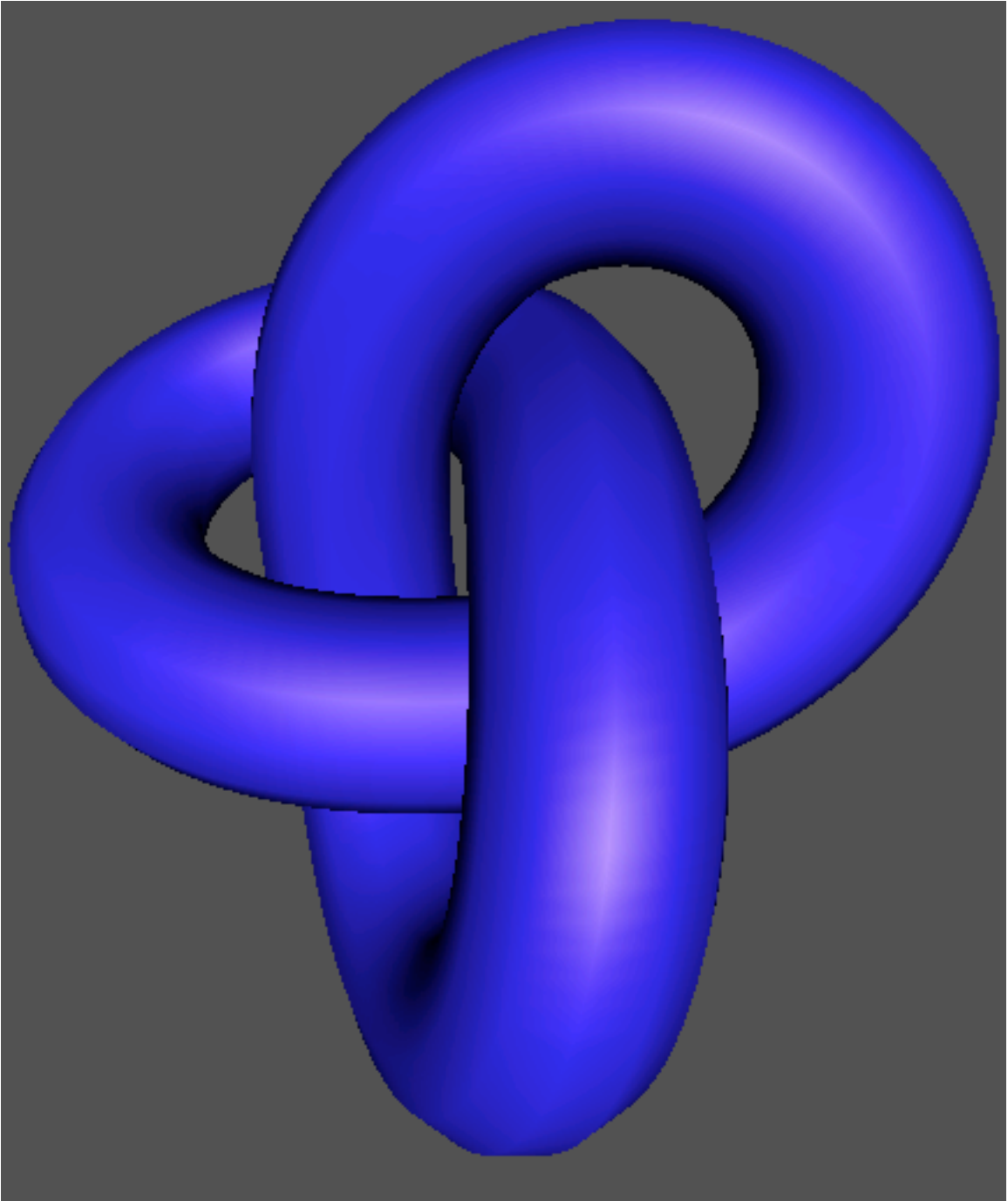

Figure 11 -- Trefoil with d = 3.5 and ERL(T, 3.5) ~ 40.45.

It is also clear that the electrical ropelength will not in general be close to the minimal ropelength for a knot. When we work with a complex knot it is often the case that the force/energy minimum for the knot has some parts of the knot considerably closer to each other than other parts. This creates an inequity that results in the electrical ropelength being greater than the minimal ropelength.

## VI. The Torus Knot (3,4) Experiment

Here is a remarkable phenomenon. When we start with a very symmetrical version of the (3,4) torus knot and use an undamped force evolution, the KnotPlot program retains the symmetry and stabilizes with the image shown in Figure 12 below. However, on perturbation via the undamped evolution, the knot falls into a lower energy level as shown in Figure 13 below. Futhermore, the electrical ropelength for d = 6 goes up when the energy falls in this example. This shows that at the experimental level, the electrical ropelength is not minimized by minimizing the Simon energy.

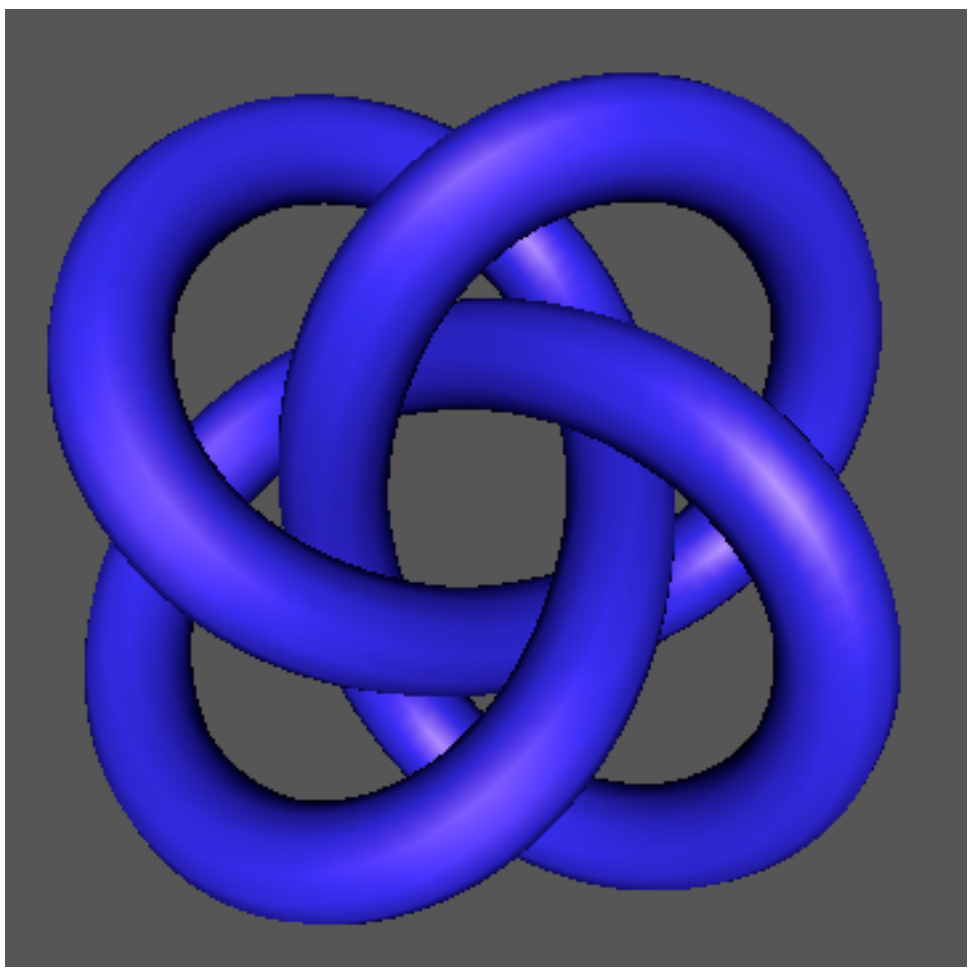

Figure 12- (3,4) Torus knot minimized via damped evolution. E = 249.306610.

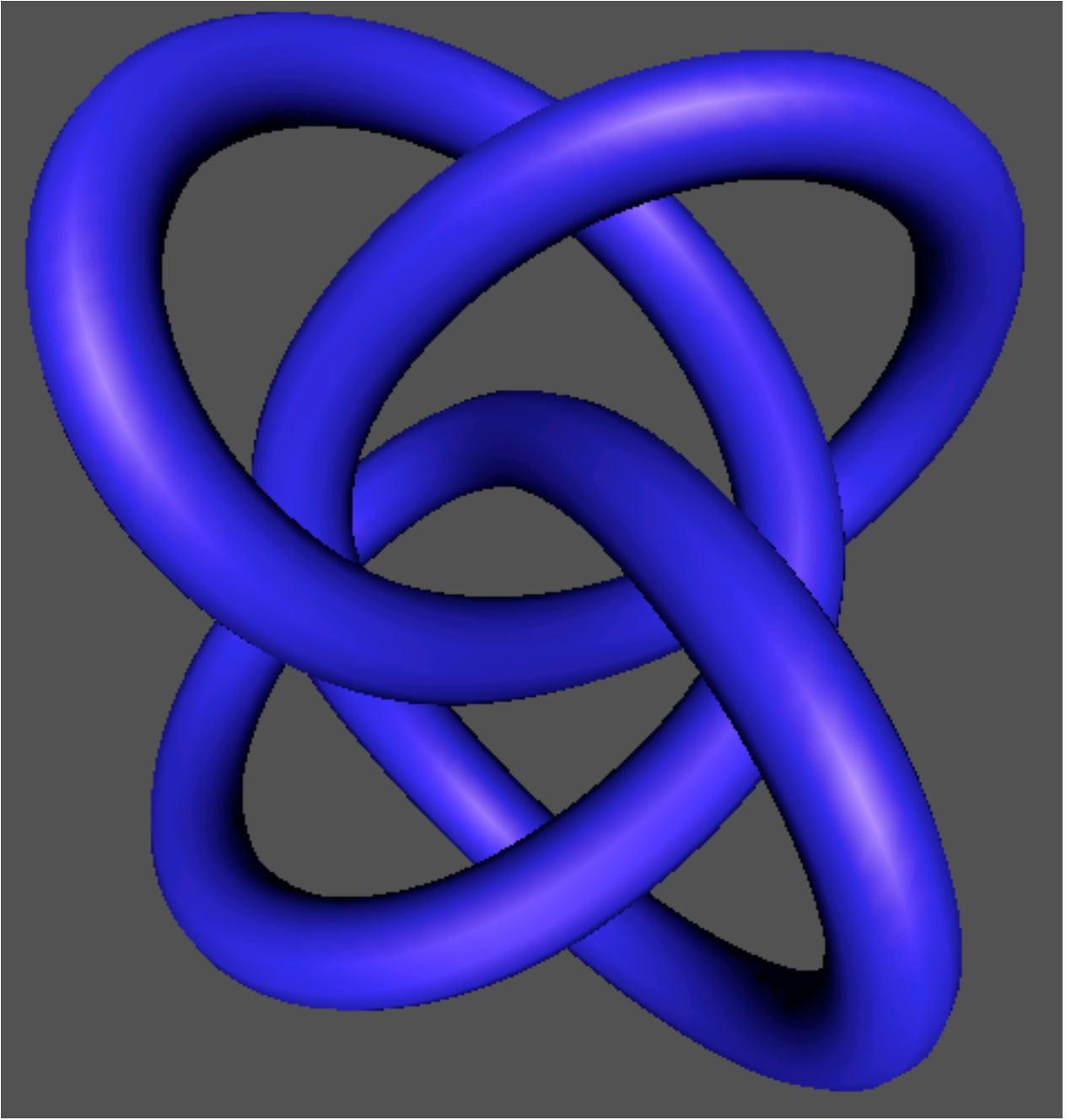

Figure 13 - (3,4) Torus knot minimized via undamped evolution. E = 247.149597.

With this example, we end out survey of experiments with energy and ropelength using the KnotPlot program.

These experiments leave many open questions. We would like to know how to better explore the landscape of configurations of a knot or link and to determine the pattern of relative and absolute minima for the Simon energy. We would like to know if, in principle, pursuing an unknot down its self-repelling evolution will unknot it (with sufficiently many beads of course). We would like to better understand the geometric configurations of the minima and their relation to their ropelengths. We encourage the reader to carry out experiments of the sort described in this paper.